\documentclass[12pt]{article}

\usepackage{amsmath, amssymb, amsfonts}

\newtheorem{theorem}{Theorem}

\newcommand{\qed}{$\Box$}

\newcommand{\op}{\operatorname}
\newcommand{\abs}[1]{\left\lvert #1 \right\rvert}

\newcommand{\br}[1]{\left\{ #1 \right\}}
\newcommand{\setdef}[2]{\br{#1 \colon #2}} 
\newcommand{\pr}[1]{\left( #1 \right)} 
\newcommand{\icard}[2]{\abs{#1 \cap #2}} 
\newcommand{\floor}[1]{\left\lfloor #1 \right\rfloor}
\newcommand{\irange}[2]{\left[#1,#2\right]} 
\newcommand{\iRange}[2]{\br{#1,\ldots,#2}} 
\newcommand{\Pb}[1]{{\mathbb{P}}{\pr{#1}}}
\newcommand{\homo}[1]{{\op{hom}}{\pr{#1}}}
\newcommand{\mono}[1]{{\op{mon}}{\pr{#1}}}

\begin{document}

\begin{center}

{\LARGE
Large cycles in generalized Johnson graphs\\
}

\vspace{0.2cm}

{\Large
Vladislav~Kozhevnikov\footnote{Moscow Institute of Physics and Technology (National Research University), Dolgoprudny, Moscow Region, Russia. Supported by Grant N NSh-775.2022.1.1 to support leading scientific schools of Russia}, Maksim~Zhukovskii\footnote{Moscow Institute of Physics and Technology (National Research University), Dolgoprudny, Moscow Region, Russia; The Russian Presidential Academy of National Economy and Public Administration, Moscow, Russia; Moscow Center for Fundamental and Applied Mathematics, Moscow, Russia;
Adyghe State University, Caucasus mathematical center, Maykop, Republic of Adygea, Russia. Supported by the Ministry of Science and Higher Education of the Russian Federation (Goszadaniye No. 075-00337-20-03), project No. 0714-2020-0005.
}
}

\vspace{0.5cm}

Abstract\\

\end{center}

We count cycles of an unbounded length in generalized Johnson graphs. Asymptotics of the number of such cycles is obtained for certain growth rates of the cycle length.



\vspace{0.5cm}

\section{Introduction and new results}
\label{intro}

For integers $i\leq j$, everywhere below we denote $\irange{i}{j}:=\iRange{i,i+1}{j}$  and $[i]:=\irange{1}{i}$. For integers $n,r,s$ such that $0\le{s}<r<n$, a simple graph $G(n,r,s)$ with the set of vertices
$$
V:=V(G(n,r,s)) = \setdef{x\subset[n]}{|x|=r}
$$
and the set of edges
$$
E:= E(G(n,r,s)) = \setdef{\br{x,y}}{\icard{x}{y} = s}
$$
is called a \textit{generalized Johnson graph}. 


Unfortunately, there is no established term for graphs $G(n,r,s)$. In literature they appear as \textit{generalized Johnson graphs} \cite{Agong2018, Cannon2012, Molitierno2017}; \textit{uniform subset graphs} \cite{Chen1987, Chen2008, Simpson1994} and \textit{distance graphs} \cite{Burkin2016, Burkin2018, Pyaderkin2016,  Zhukovskii2012_sub}. The family of $G(n,r,s)$ graphs was initially (to the best of our knowledge) considered in \cite{Chen1987}, where they are called ``\textit{uniform subset graphs}''. However, this name does not become widespread. In our opinion, the term ``\textit{generalized Johnson graph}'' is preferred as the most comprehensible, since, if we set $s=r-1$, then the definition of $G(n,r,s)$ turns into the definition of the well-known Johnson graph. Note that the Kneser graph is also a special case of $G(n,r,s)$ with $s=0$. However, the term ``\textit{generalized Kneser graph}" is already used for another class of graphs \cite{Chen2008_generalized_kneser, Denley1997, Frankl1985_generalized_kneser,  Jafari2020}.



On the one hand, as we mentioned above, graphs $G(n,r,s)$ generalize Johnson graphs $G(n,r,r-1)$ \cite{Alspach2012, Daven1999, Etzion1996_codes, Etzion1996_chromatic, Etzion2011} and Kneser graphs $G(n,r,0)$ \cite{Brouwer1979, Chen2003, Lovasz1978, Matousek2004, Mutze2020, Poljak1987, Valencia2005}, which are themselves of interest in the graph theory. On the other hand, they are a special case of distance graphs in $\mathbb{R}^n$ with the Euclidean metric, which are used to study problems of combinatorial geometry (Hadwiger--Nelson problem about the chromatic number $\chi(\mathbb{R}^n)$ \cite{Cantwell1996,  Chilakamarri1990, Exoo2014, Frankl1980, Frankl1981, Kupavskii2009, Larman1972}, Borsuk problem about partitioning of a set in $\mathbb{R}^n$ into subsets of a smaller diameter \cite{Kahn1993}, and various generalizations of these problems \cite{Raigorodskii2001, Raigorodskii2016,  Raigorodskii2013}).\\ 


Throughout the paper we assume that $r$ and $s$ are constant and $n$ approaches  infinity. The total number of vertices in this graph is denoted by $N$:
\[N = \abs{V} = \binom{n}{r} \sim \frac{n^r}{r!}.\]
From the definition of $G(n,r,s)$ it is evident that this graph is \textit{vertex--transitive}, i.e. for any two vertices there exists an automorphism of the graph mapping the first vertex to the second one. In particular, $G(n,r,s)$ is regular. Let $N_1$ denote the degree of its vertex: \[N_1 = \binom{r}{s} \binom{n-r}{r-s} \sim \binom{r}{s} \frac{n^{r-s}}{(r-s)!}.\]

In \cite{Chen1987} it is proved that the graph $G(n,r,s)$ is Hamiltonian for $s\in\br{r-1,r-2,r-3}$, arbitrary $r$ and $n$ as well as for $s\in\br{0,1}$, arbitrary $r$, and sufficiently large $n$.
Hamiltonian cycles has been extensively studied in Kneser graphs $G(n,r,0)$. It is known that they are hamiltonian for $n\ge2.62r$ \cite{Chen2003} and for all $r$ when $n\le27$ (except for the Petersen graph $G(5,2,0)$) \cite{Shields2002}.
Graphs $G(2r+1,r,0)$ are also known to be Hamiltonian for all $r\ge3$ \cite{Mutze2020}.
As for cycles of a constant length, the asymptotics of the number of their appearances in $G(n,r,s)$ is known for all constant $r$ and $s$ and given below in Theorem~\ref{th_fixed_t_mon_asymp}.


Let $H$ and $G$ be graphs. A map $\varphi:V(H)\to{V(G)}$ is called a \textit{homomorphism} from $H$ to $G$ if, for any pair of vertices $x,y$ of $H$, $\br{x,y}\in{E(H)}\Rightarrow\br{\varphi(x),\varphi(y)}\in{E(G)}$. If a homomorphism is injective, then it is called a \textit{monomorphism}. Let $\homo{H,G}$ and $\mono{H,G}$ denote respectively the number of homomorphisms and monomorphisms from $H$ to $G$. Throughout this paper we write simply $\homo{H}$ and $\mono{H}$ when $G=G(n,r,s)$.

Let $C_t$ be a cycle on $t$ vertices.
The purpose of this paper is to find the asymptotic value of $\mono{C_t}$ for different $t=t(n)$.

Burkin \cite{Burkin2016} found the asymptotics of $\mono{C_t}$ for all $t=\op{const}$.

\begin{theorem}[Burkin, 2016, \cite{Burkin2016}]
\label{th_fixed_t_mon_asymp}
Let $t$ be a fixed integer. Then
\begin{equation}
\label{eq_fixed_t_mon_asymp}
\mono{C_t} \sim N N_1 \pr{\frac{N_1}{\binom{r}{s}}}^{t-2}.
\end{equation}
\end{theorem}
We generalize this result to cycles of variable length, i.e. $t=t(n)$. It turns out that for slow enough (sublogarithmic) growth of $t(n)$ the asymptotics of $\mono{C_t}$ remains the same as in \eqref{eq_fixed_t_mon_asymp}. In contrast, for superlogarithmic $t(n)=o\pr{\min\{\sqrt{N},N_1\}}$ the asymptotics is different, namely, $\mono{C_t}\sim N_1^t$. These results can be summarized in the following two theorems.

\begin{theorem}
\label{th_mon_asymp_eq_hom}
As $n\to+\infty$,
$\mono{C_t}\sim\homo{C_t}$
iff $t=o\pr{\min\{\sqrt{N},N_1\}}$.
\end{theorem}

Theorem~\ref{th_mon_asymp_eq_hom} is the trickiest result of our paper. Asymptotics of $\homo{C_t}$ (stated below in Theorem~\ref{th_hom}) is a more or less direct corollary (modulo technical asymptotical computations) of the well-known representation of $\homo{C_t}$ in terms of eigenvalues of $G(n,r,s)$. Let us fix an arbitrarily small $\varepsilon>0$ and consider the partition of $\mathbb{N}$ obtained by excluding $\varepsilon$-neighborhoods of $\frac{\ln n}{\ln{r-j\choose s-j}-\ln{r-j-1\choose s-j-1}}$, $j\in[0,s]$, i.e. the intervals 
$$
I_j=\left[\left\lfloor\frac{(1+\varepsilon)\ln{n}}{\ln{r-j\choose s-j}-\ln{r-j-1\choose s-j-1}}\right\rfloor\right]\setminus\left[\left\lfloor\frac{(1-\varepsilon)\ln{n}}{\ln{r-j+1\choose s-j+1}-\ln{r-j\choose s-j}}\right\rfloor\right],\quad j\in[s-1],
$$
$$
I_s=\left[\left\lfloor\frac{(1-\varepsilon)\ln{n}}{\ln(r-s+1)}\right\rfloor\right],\quad I_0=\left[\left\lfloor\frac{(1+\varepsilon)\ln{n}}{\ln\frac{r}{s}}\right\rfloor,\infty\right).
$$

\begin{theorem}
\label{th_hom}
For arbitrary $t=t(n)\in\mathbb{N}$,
\begin{equation}
\label{eq_hom}
\homo{C_t}={N_1^t}\pr{1+O\left(\frac{1}{n}\right)+\sum\limits_{j=1}^{s}\frac{n^j}{j!}\pr{\frac{\binom{r-j}{s-j}}{\binom{r}{s}}+O\left(\frac{1}{n}\right)}^t}.    
\end{equation}
Moreover, for $j\in[0,s+1]$ and $t\in I_j$,
$$
\homo{C_t}\sim N_1^t\frac{n^j}{j!}\left({r-j\choose s-j}/{r\choose s}\right)^t.
$$
\end{theorem}

Note that for $t\in I_0$, $\homo{C_t}\sim N_1^t$, while, for $t\in I_s$,
$\homo{C_t}\sim N_1^t\frac{n^s}{s!}{r\choose s}^{-t}$
i.e. \eqref{eq_fixed_t_mon_asymp} holds. Theorem~\ref{th_mon_asymp_eq_hom} and Theorem~\ref{th_hom} immediately yield asymptotics of the number of copies of $C_t$ in $G(n,r,s)$ for all $t=o\pr{\min\{\sqrt{N},N_1\}}$ since it equals $\frac{1}{2t}\mono{C_t}$.\\



The rest of the paper is organized as follows. First, in Section~\ref{random_walk}, we discuss general properties of random walks on graphs (Section~\ref{sec:trans_matrix_and_mixing_time}) and more specific properties of random walks on $G(n,r,s)$ (Sections~\ref{sec:eigen_Johnson} and \ref{sec:rw_Janson}).
Secondly, in Section~\ref{proof_hom}, we prove Theorem~\ref{th_hom}. Finally, in Sections~\ref{proof_mon_asymp_eq_hom} and \ref{sec:th3_proof} we prove that the condition in Theorem~\ref{th_mon_asymp_eq_hom} is, respectively, sufficient and necessary.

The proof of the sufficiency provided in Section~\ref{proof_mon_asymp_eq_hom} uses exact expressions for the spectrum of $G(n,r,s)$. It should be noted that the proof in the case $r>2s$ (in which $\sqrt{N}=o(N_1)$) as well as in the case $t=\omega(\ln{N})$ can be considerably simplified by using a more general argument (which we omit in this paper) applicable to a wide subclass of spectral expanders (see Section~\ref{sec:th3_proof}).
However, for an arbitrary $N_1$-regular graph $G$ on $N$ vertices, the property $\mono{C_t,G}\sim\homo{C_t,G}$ does not necessarily hold when
$t=O(\ln{N})$
and $N_1=O(\sqrt{N})$, even if $G$ is a spectral expander.
This fact can be demonstrated, for example, by considering the random regular graph $G(N,N_1)$ with $N_1=\lfloor\ln^8{N}\rfloor$, in which, for any $\varepsilon>0$, the inequality $\mono{C_t,G(N,N_1)}/\homo{C_t,G(N,N_1)}<\varepsilon$ holds with probability approaching $1$ as soon as $t=o(\ln{N}/\ln\ln{N})$. This can be shown by translating the same property from the binomial random graph $G(N,(1+o(1))N_1/N)$ to $G(N,N_1)$ using the sandwich conjecture, which is true for $N_1=\omega(\ln^7{N})$~\cite{Gao2021}.
Note that $G(N,N_1)$ is a spectral expander \cite{Zhao2012}.
For the definition and properties of binomial random graphs and regular random graphs see \cite{Janson2000}.

The proof of the necessity in Theorem~\ref{th_mon_asymp_eq_hom} provided in Section~\ref{sec:th3_proof} does not rely upon the the whole spectrum of $G(n,r,s)$ but rather uses its spectral expansion property.
The necessity of the condition $t=o(N_1)$ follows from the fact that a random walk starts backtracking with positive probability if $t>cN_1$ for a constant $c$, which is proved in Section~\ref{sec:th3_proof_2} using almost solely the regularity of $G(n,r,s)$.
The necessity of $t=o(\sqrt{N})$ is proved in Section~\ref{sec:th3_proof_1} using a high convergence rate of a random walk on an expander, which is discussed in Section~\ref{sec:trans_matrix_and_mixing_time}. Therefore, in Section~\ref{sec:th3_proof} we formulate a generalization of Theorem~\ref{th_mon_asymp_eq_hom} to a class of spectral expanders.

\section{Random walks}
\label{random_walk}

Counting cycles in $G(n,r,s)$ can be reduced to analysing the distribution of a random walk on $G(n,r,s)$.

\subsection{Distribution and adjacency matrix}
\label{sec:trans_matrix_and_mixing_time}

Let $G$ be an arbitrary regular connected graph on the vertex set $[N]$ with every vertex having degree $N_1$. Let $A=(A_{i,j},\,i,j\in[N])$ be its adjacency matrix ($A_{ij}=1$ if and only if $i$ and $j$ are adjacent in $G$). Moreover, let $\lambda_j$, $j\in[0,r]$, be all distinct eigenvalues of $A$, and let $m_j$ be the multiplicity of $\lambda_j$.\\

Recall that {\it a random walk on $G$} is a discrete-time random process $(X_n,\,n\in\mathbb{Z}_+)$, where $X_0$ is a vertex chosen uniformly at random from $[N]$, and, for every $n\in\mathbb{Z}_+$, $X_{n+1}$ is chosen uniformly at random from the neighbors of $X_n$ in $G$. For $x,y\in[N]$, let 
$$
P^t(x,y):=\Pb{X_t=y|X_0=x}
$$
and $P^t=(P^t(x,y),\,x,y\in[N])$ be the $k$-step transition probabilities matrix.


For a positive integer $t$, $\homo{C_t,G}$ is exactly the trace of $A^t$. Since the trace of $A^t$ equals the sum of its eigenvalues and the eigenvalues of $A^t$ can be computed as the $t$th power of the eigenvalues of $A$~(see, e.g., \cite{Meyer2000}), we get 
\begin{equation}
\homo{C_t,G}=\sum\limits_{j=0}^{r}m_j\lambda_j^t.
\label{eq:homo_eigen}
\end{equation}

If $G$ is vertex-transitive, then, clearly, all $P^t(z,z)$, $z\in[N]$, are equal to each other. Then, for every $x\in[N]$,
\begin{equation}
NP^t(x,x)=\sum_{z\in[N]} P^t(z,z)=\frac{\homo{C_t,G}}{N_1^t}.
\label{eq:homo_prob}
\end{equation} 
Therefore,~(\ref{eq:homo_eigen}) implies that
\begin{equation}
\label{eq_t_step_trans_prob}
P^t(x,x)=\frac{1}{N}\sum\limits_{j=0}^{r}m_j\pr{\frac{\lambda_j}{N_1}}^t.    
\end{equation}

Notice that (due to regularity of $G$) the distribution $\pi=(1/N,\ldots,1/N)\in\mathbb{R}^N$ is {\it stationary} meaning that $\pi P^1=\pi$. Let us here assume that $\lambda_0$ is the largest eigenvalue and $\lambda_1$ is the largest in absolute value eigenvalue distinct from $\lambda_0$.
From the regularity of $G$ it follows that $\lambda_0=N_1$ and and from its connectedness, that $m_0=1$~\cite{Bapat2014}.
Let us also assume that $\abs{\lambda_1}<\lambda_0$ (which is equivalent, for a connected graph, to the graph being non-bipartite~\cite{Bapat2014}).
Fix $v\in [N]$ and $\varepsilon>0$. Let us recall that {\it the variation distance} at time $t\in\mathbb{Z}_+$ with initial state $v$ is
$$
 \Delta_v(t)=\frac{1}{2}\sum_{u\in[N]}\left|P^t(v,u)-\frac{1}{N}\right|.
$$
It is very well known~\cite{Sinclair1992} that {\it the mixing time} $\tau_v(\varepsilon):=\min\{t:\,\Delta_v(t')\leq\varepsilon\text{ for all }t'\geq t\}$ satisfies
\begin{equation}
\tau_v(\varepsilon)\leq\left(1-\frac{|\lambda_1|}{N_1}\right)^{-1}\ln\frac{N}{\varepsilon}.
\label{eq:mix_time}
\end{equation}

\subsection{Eigenvalues of $G(n,r,s)$}
\label{sec:eigen_Johnson}


Let $A$ be the adjacency matrix of $G(n,r,s)$. The eigenvalues of $A$ are known~\cite{Delsarte1973}. They are equal to
\begin{equation}
\lambda_j=\sum\limits_{\ell=\max\{0,j-s\}}^{\min\{j,r-s\}}(-1)^{\ell}\binom{j}{\ell}\binom{r-j}{r-s-\ell}\binom{n-r-j}{r-s-\ell}, \quad j\in[0,r],
\label{eq:eigen_dist_exact}
\end{equation}
and the multiplicity of the eigenvalue $\lambda_j$ equals (we let ${n\choose -1}=0$)
$$
m_j=\binom{n}{j}-\binom{n}{j-1}.
$$

In order to prove Theorem~\ref{th_hom}, we need to analyze asymptotical behaviour of the expression to the right in (\ref{eq_t_step_trans_prob}).

Notice that $\lambda_0=N_1$ and $m_0=1$. Also, for $j\in\irange{1}{s}$,
\begin{equation}
\frac{\lambda_j}{N_1}=\frac{\binom{r-j}{s-j}}{\binom{r}{s}}+O\left(\frac{1}{n}\right),\quad m_j=\frac{n^j}{j!}+O\left(\frac{1}{n}\right)
\label{eq:eigen_distance_asymp1}
\end{equation}
and, for $j\in\irange{s+1}{r}$,
\begin{equation}
\abs{\frac{\lambda_j}{N_1}}\sim\frac{\binom{j}{s}(r-s)!}{\binom{r}{s}(r-j)!}n^{-(j-s)},\quad m_j=\frac{n^j}{j!}+O\left(\frac{1}{n}\right).
\label{eq:eigen_distance_asymp2}
\end{equation}

Therefore, for $t\ge2$ and $j\in\irange{s+1}{r}$,
\begin{equation}
\frac{m_j\lambda_j^t}{m_s\lambda_s^t}=(O(1))^t\cdot{n^{-(j-s)(t-1)}}=\pr{O\pr{n^{-(j-s)}}}^{t-1}=O\pr{\frac{1}{n}}
\label{eq:compare_eigen_dist_large_j}
\end{equation}
implying that
\begin{equation}
P^t(x,x)=\frac{1+O\pr{\frac{1}{n}}}{N}\sum\limits_{j=0}^{s}m_j\pr{\frac{\lambda_j}{N_1}}^t=
\frac{1}{N}\left(1+O\left(\frac{1}{n}\right)+\sum\limits_{j=1}^{s}m_j\left(\frac{\binom{r-j}{s-j}}{\binom{r}{s}}+O\left(\frac{1}{n}\right)\right)^t\right).    
\label{eq:hom_eigen_asymp}
\end{equation}

\subsection{Random walk on $G(n,r,s)$}
\label{sec:rw_Janson}

Here, we consider the random walk $(X_n,\,n\in\mathbb{Z}_+)$ on $G(n,r,s)$. Since $G(n,r,s)$ is vertex-transitive, for any $x\in{V}$, 
$$
\mono{C_t}=NN_1^t\Pb{X_t=x,X_0\ne{X_1}\ne\ldots\ne{X_{t-1}}|X_0=x}.
$$
In order to prove Theorem~\ref{th_mon_asymp_eq_hom}, we bound the deviation of $\frac{\mono{C_t}}{\homo{C_t}}$ from $1$ . For convenience, in this section, we express the bound in terms of the diagonal elements of $P^t$.

Due to (\ref{eq:homo_prob}), we get 

\begin{align*}
0\le\frac{\homo{C_t}-\mono{C_t}}{\homo{C_t}}
&=\frac{P^t(x,x)-\Pb{X_t=x,X_0\ne{X_1}\ne\ldots\ne{X_{t-1}}|X_0=x}}{P^t(x,x)}=\\
&=\frac{\Pb{X_t=x,\exists{i,j}\in\irange{0}{t-1}:i\ne{j},X_i=X_j|X_0=x}}{\Pb{X_t=x|X_0=x}}.
\end{align*}
Note that the expression to the right is exactly the probability that the random walk meets itself somewhere on $[0,t-1]$ subject to $X_0=x$ and $X_t=x$.

By the union bound,
\begin{align*}
\frac{\homo{C_t}-\mono{C_t}}{\homo{C_t}}&\le\sum\limits_{0\le{i}<j<t}\frac{\Pb{X_j=X_i,X_t=x|X_0=x}}{\Pb{X_t=x|X_0=x}}\\
&=\sum\limits_{0\le{i}<j<t}\frac{\sum\limits_{z\in{V}}\Pb{X_t=x|X_j=z}\Pb{X_j=z|X_i=z}\Pb{X_i=z|X_0=x}}{\Pb{X_t=x|X_0=x}}.
\end{align*}

Due to vertex-transitivity of $G(n,r,s)$ the probabilities $\Pb{X_j=z|X_i=z}$ are equal for all $z$. Therefore,
\begin{equation}
\begin{split}
\frac{\homo{C_t}-\mono{C_t}}{\homo{C_t}}
&\leq\sum\limits_{0\le{i}<j<t}\frac{\sum\limits_{z\in{V}}\Pb{X_{t-j+i}=x|X_i=z}P^{j-i}(x,x)\Pb{X_i=z|X_0=x}}{\Pb{X_t=x|X_0=x}}\\
&=\sum\limits_{0\le{i}<j<t}\frac{P^{t-j+i}(x,x)P^{j-i}(x,x)}{P^{t}(x,x)}
=\sum\limits_{k=1}^{t-1}(t-k)\frac{P^{k}(x,x)P^{t-k}(x,x)}{P^{t}(x,x)}\\
&\le{t}\sum\limits_{k=1}^{t-1}\frac{P^{k}(x,x)P^{t-k}(x,x)}{P^{t}(x,x)}={t}\sum\limits_{k=2}^{t-2}\frac{P^{k}(x,x)P^{t-k}(x,x)}{P^{t}(x,x)}.
\label{eq_gap_upper_bound}
\end{split}
\end{equation}

\section{Proof of Theorem~\ref{th_hom}}
\label{proof_hom}
From (\ref{eq:homo_prob}) and (\ref{eq:hom_eigen_asymp}), we get
\begin{equation}
\homo{C_t}=
N_1^t\pr{1+O\left(\frac{1}{n}\right)+\sum\limits_{j=1}^{s}\frac{n^j}{j!}\pr{\frac{\binom{r-j}{s-j}}{\binom{r}{s}}+O\left(\frac{1}{n}\right)}^t}.
\label{eq:hom_expansion}
\end{equation}
Let $j\in[s-1]$, $\varepsilon>0$. If $t<(1-\varepsilon)\frac{\ln n}{\ln{r-j\choose s-j}-\ln{r-j-1\choose s-j-1}}$, then
$$
 \left(\frac{{r-j\choose s-j}}{{r\choose s}}+o(1)\right)^t=
 o\left[n\left(\frac{{r-j-1\choose s-j-1}}{{r\choose s}}+o(1)\right)^t\right]
$$
implying that the $(j+1)$th term in the sum in (\ref{eq:hom_expansion}) is asymptotically bigger than the $j$th term. Similarly, if $t>(1+\varepsilon)\frac{\ln n}{\ln{r-j\choose s-j}-\ln{r-j-1\choose s-j-1}}$, then the $(j+1)$th term in the sum in (\ref{eq:hom_expansion}) is asymptotically smaller than the $j$th term. From this, 
we immediately get

$$
\homo{C_t}\sim
\begin{cases}
N_1^t\frac{n^s}{s!}{r\choose s}^{-t}, & t<(1-\varepsilon)\frac{\ln n}{\ln{r-s+1\choose s-s+1}-\ln{r-s\choose s-s}};  \\
N_1^t\frac{n^j}{j!}{r-j\choose s-j}^t{r\choose s}^{-t}, & \frac{(1+\varepsilon)\ln n}{\ln{r-j\choose s-j}-\ln{r-j-1\choose s-j-1}}<t<\frac{(1-\varepsilon)\ln n}{\ln{r-j+1\choose s-j+1}-\ln{r-j\choose s-j}},\,j\in[s-1];  \\
N_1^t, & t>(1+\varepsilon)\frac{\ln n}{\ln{r\choose s}-\ln{r-1\choose s-1}}.
\end{cases}
$$
which proves Theorem~\ref{th_hom}.
\qed

\section{Proof of sufficiency in Theorem~\ref{th_mon_asymp_eq_hom}}
\label{proof_mon_asymp_eq_hom}

Here we prove that if $t=o\pr{\min\{\sqrt{N},N_1\}}$, then
$$
\frac{\homo{C_t}-\mono{C_t}}{\homo{C_t}}\to 0,\quad n\to\infty.
$$
Since this fraction is non-negative, it is sufficient to prove that the upper bound from (\ref{eq_gap_upper_bound}) approaches 0. Clearly, we may assume that $t\geq 4$.\\

By (\ref{eq_t_step_trans_prob}), for every $x\in V$,
\begin{equation}
\begin{split}
\sum\limits_{k=2}^{t-2}P^{k}(x,x)P^{t-k}(x,x)
&=\frac{1}{N^2}\sum\limits_{k=2}^{t-2}\sum\limits_{i,j=0}^{r}m_im_j
\left(\frac{\lambda_j}{N_1}\right)^k\left(\frac{\lambda_i}{N_1}\right)^{t-k}\\
&\le\frac{t}{N^2}\sum\limits_{i=0}^{r}m_i^2\abs{\frac{\lambda_i}{N_1}}^t+
\frac{1}{N^2}\sum\limits_{i\ne{j}}{m_i{m_j}\abs{\frac{\lambda_i}{N_1}}^t\sum\limits_{k=2}^{t-2}\abs{\frac{\lambda_j}{\lambda_i}}^{k}}.
\end{split}
\label{eq:diagonal_decompose}
\end{equation}
Note that $\abs{\lambda_i}^t\sum\limits_{k=2}^{t-2}\abs{\frac{\lambda_j}{\lambda_i}}^{k}=\abs{\lambda_j}^t\sum\limits_{k=2}^{t-2}\abs{\frac{\lambda_i}{\lambda_j}}^{k}$. Therefore, for $n$ large enough, due to (\ref{eq:eigen_distance_asymp1})~and~(\ref{eq:eigen_distance_asymp2}), we get
$$
\sum\limits_{i\ne{j}}{m_i{m_j}\abs{\frac{\lambda_i}{N_1}}^t\sum\limits_{k=2}^{t-2}\abs{\frac{\lambda_j}{\lambda_i}}^{k}}=
2\sum\limits_{0\le{i}<j\le{r}}{m_i{m_j}\abs{\frac{\lambda_i}{N_1}}^t\sum\limits_{k=2}^{t-2}\abs{\frac{\lambda_j}{\lambda_i}}^{k}}\leq
2\sum\limits_{0\le{i}<j\le{r}}{m_i{m_j}\abs{\frac{\lambda_i}{N_1}}^t\frac{\abs{\frac{\lambda_j}{\lambda_i}}^{2}}{1-\abs{\frac{\lambda_j}{\lambda_i}}}}.
$$
Moreover, from (\ref{eq:eigen_distance_asymp1})~and~(\ref{eq:eigen_distance_asymp2}) we get that, for $j>i$ and $j>s$, $\frac{\lambda_j}{\lambda_i}=O\left(\frac{1}{n}\right)$, while, for $i<j\leq s$, 
$$
\frac{\lambda_j}{\lambda_i}=\frac{(s-i)\ldots(s-j+1)}{(r-i)\ldots(r-j+1)}\left(1+O\left(\frac{1}{n}\right)\right).
$$
The latter expression is less than $\frac{s}{r}$ if $j\neq 1$ and $n$ is large enough. If $i=0$, $j=1$, then, from (\ref{eq:eigen_dist_exact}), we get
$$
\frac{\lambda_1}{\lambda_0}=\frac{{r-1\choose r-s}{n-r-1\choose r-s}-{r-1\choose r-s-1}{n-r-1\choose r-s-1}}{{r\choose r-s}{n-r\choose r-s}}<
\frac{{r-1\choose r-s}{n-r-1\choose r-s}}{{r\choose r-s}{n-r\choose r-s}}=\frac{s(n-2r+s)}{r(n-r)}<\frac{s}{r}.
$$
Then, for $n$ large enough,
\begin{align*}
\sum\limits_{i\ne{j}}{m_i{m_j}\abs{\frac{\lambda_i}{N_1}}^t\sum\limits_{k=2}^{t-2}\abs{\frac{\lambda_j}{\lambda_i}}^{k}}
&\leq\frac{2}{1-s/r}\sum\limits_{0\le{i}<j\le{r}}{m_i{m_j}\abs{\frac{\lambda_i}{N_1}}^{t-2}\abs{\frac{\lambda_j}{N_1}}^{2}}\\
&\leq\frac{2}{1-s/r}\sum\limits_{i,j=0}^r{m_i{m_j}\abs{\frac{\lambda_i}{N_1}}^{t-2}\abs{\frac{\lambda_j}{N_1}}^{2}}\\
&=\frac{2}{1-s/r}\left(\sum\limits_{i=0}^r m_i\abs{\frac{\lambda_i}{N_1}}^{t-2}\right)\left(\sum\limits_{j=0}^r m_j\abs{\frac{\lambda_j}{N_1}}^{2}\right).
\end{align*}
By (\ref{eq:hom_eigen_asymp}) and the definition of $P^2$, 
$$
\sum\limits_{j=0}^r m_j\abs{\frac{\lambda_j}{N_1}}^{2}=NP^2(x,x)=\frac{N}{N_1}.
$$
Moreover, by (\ref{eq:eigen_distance_asymp1}), (\ref{eq:compare_eigen_dist_large_j}) and (\ref{eq:hom_eigen_asymp}),
\begin{align*}
 \sum\limits_{i=0}^r m_i\abs{\frac{\lambda_i}{N_1}}^{t-2}\sim
 \sum\limits_{i=0}^s m_i\abs{\frac{\lambda_i}{N_1}}^{t-2}\sim
 \sum\limits_{i=0}^s m_i\left(\frac{\lambda_i}{N_1}\right)^{t-2}
 &=O\left(\sum\limits_{i=0}^sm_i\left(\frac{\lambda_i}{N_1}\right)^{t}\right)\\
 &=O\left(\sum\limits_{i=0}^r m_i\left(\frac{\lambda_i}{N_1}\right)^t\right)=
 O(NP^t(x,x)).
\end{align*}
It remains to estimate the first summand in the rightmost expression in (\ref{eq:diagonal_decompose}). From (\ref{eq:eigen_distance_asymp1})~and~(\ref{eq:eigen_distance_asymp2}), we get that, for $j\in\irange{s+1}{r}$,
$$
\frac{m_j^2\lambda_j^t}{m_s^2\lambda_s^t}=(O(1))^t\cdot{n^{-(j-s)(t-2)}}=\pr{O\pr{n^{-(j-s)}}}^{t-2}=o(1)
$$
implying that
$$
\sum\limits_{i=0}^{r}m_i^2\abs{\frac{\lambda_i}{N_1}}^t\sim
\sum\limits_{i=0}^{s}m_i^2\left(\frac{\lambda_i}{N_1}\right)^t.
$$
Putting everything together and applying (\ref{eq:hom_eigen_asymp}), we conclude that, for every $x\in V$,
$$
t\sum\limits_{k=2}^{t-2}\frac{P^{k}(x,x)P^{t-k}(x,x)}{P^{t}(x,x)}\leq
\frac{t^2}{N}\frac{\sum\limits_{i=0}^{s}m_i^2\left(\frac{\lambda_i}{N_1}\right)^t}{\sum\limits_{i=0}^{s}m_i \left(\frac{\lambda_i}{N_1}\right)^t}(1+o(1))+
O\left(\frac{t}{N_1}\right).
$$
We finish with proving that the condition $t=o(\sqrt{N})$ implies
$$
 \frac{t^2}{N}\frac{\sum\limits_{i=0}^{s}m_i^2\left(\frac{\lambda_i}{N_1}\right)^t}{\sum\limits_{i=0}^{s}m_i \left(\frac{\lambda_i}{N_1}\right)^t}=o(1).
$$
If $t>C\ln{n}$, where $C$ is sufficiently large constant, then, by (\ref{eq:eigen_distance_asymp1}), for every $i\in[s]$, $m_i\left(\frac{\lambda_i}{N_1}\right)^t=o(1)$ and $m_i^2\left(\frac{\lambda_i}{N_1}\right)^t=o(1)$. Therefore,
$$
 \frac{t^2}{N}\frac{\sum\limits_{i=0}^{s}m_i^2\left(\frac{\lambda_i}{N_1}\right)^t}{\sum\limits_{i=0}^{s}m_i \left(\frac{\lambda_i}{N_1}\right)^t}\sim\frac{t^2}{N}=o(1).
$$
If $t\leq C\ln{n}$, then, by (\ref{eq:eigen_distance_asymp1}),
$$
 \frac{t^2}{N}\frac{\sum\limits_{i=0}^{s}m_i^2\left(\frac{\lambda_i}{N_1}\right)^t}{\sum\limits_{i=0}^{s}m_i \left(\frac{\lambda_i}{N_1}\right)^t}\leq\frac{t^2 m_s}{N}(1+o(1))=O\left(\frac{n^s\ln^2 n}{N}\right)=O\left(\frac{\ln^2 n}{n^{r-s}}\right)=o(1).
$$
The sufficiency in Theorem~\ref{th_mon_asymp_eq_hom} is proved.
\qed

\section{Proof of necessity in Theorem~\ref{th_mon_asymp_eq_hom}}
\label{sec:th3_proof}

As was already noted, the necessity of the condition in Theorem~\ref{th_mon_asymp_eq_hom} follows from a more general fact about spectral expanders. Consider a sequence of graphs $\{G_N, N\in\mathbb{N}\}$ such that $[N]$ is the set of vertices of $G_N$, and $G_N$ is non-bipartite, connected and $N_1$-regular ($N_1$ depends on $N$). Let $\lambda_1=\lambda_1(N)$ be the second largest in absolute value eigenvalue of the adjacency matrix of $G_N$. We call the sequence $\{G_N, N\in\mathbb{N}\}$ a
\textit{spectral expander}, if there exists $\delta>0$ such that, for all large enough $N$, $\abs{\lambda_1}/N_1<1-\delta$.

\begin{theorem}
\label{th_mon_asymp_non_eq_hom}
Let $\{G_N, N\in\mathbb{N}\}$ be a spectral expander such that $N_1=\omega(\ln{N})$. Then for any $c>0$ there exists $\varepsilon>0$ such that, for sufficiently large $N$, if $t>c\min\{\sqrt{N},N_1\}$, then $\mono{C_t,G_N}<(1-\varepsilon)\homo{C_t,G_N}$.
\end{theorem}

Obviously, Theorem~\ref{th_mon_asymp_non_eq_hom}
implies the necessity of the condition $t=o\pr{\min\{\sqrt{N},N_1\}}$ for $\mono{C_t}\sim\homo{C_t}$ stated in Theorem~\ref{th_mon_asymp_eq_hom}.\\

To prove Theorem~\ref{th_mon_asymp_non_eq_hom}, we introduce a random walk on $G_N$ (see Section~\ref{sec:trans_matrix_and_mixing_time}).

In Section~\ref{sec:rw_Janson}, we note that the proportion
of self-intersecting cycles in $G(n,r,s)$ is exactly the probability that the random walk meets itself somewhere on $[0,t-1]$ subject to $X_0=x$ and $X_t=x$. It is easy to see that (almost) the same is true for $G_N$ since $X_0$ is chosen uniformly at random from $[N]$. Indeed,
\begin{align*}
    \homo{C_t,G_N}=\sum_{x\in[N]}P^t(x,x)N_1^t
    &=NN_1^t\sum_{x\in[N]}P^t(x,x)\Pb{X_0=x}\\
    &=NN_1^t\sum_{x\in[N]}\Pb{X_t=X_0=x}=NN_1^t\Pb{X_t=X_0}.
\end{align*}
In the same way, $\mono{C_t,G}=NN_1^t\Pb{X_t=X_0,X_0\neq X_1\neq\ldots\neq X_{t-1}}$. Therefore,
\begin{equation}
    \frac{\homo{C_t,G}-\mono{C_t,G}}{\homo{C_t,G}}=\Pb{\exists{i,j}\in\irange{0}{t-1}:i\ne{j},X_i=X_j|X_t=X_0}.
\label{eq:relative_error}
\end{equation}

Further, we separately consider cases $\sqrt{N}=o(N_1)$
and $N_1=O(\sqrt{N})$.

\subsection{$\sqrt{N}=o(N_1)$}
\label{sec:th3_proof_1}

Let us first bound from above $P^t(x,y)$ for arbitrary $x,y$ from $[N]$ and an integer $t\geq 1$. Since $P^t(x,y)=\sum_{v\in [N]}P^{t-1}(x,v)P^1(v,y)$, we get that 
\begin{equation}
P^t(x,y)\leq\max_{v\in[N]} P^1(v,y)\leq\frac{1}{N_1}.
\label{eq:transition_above}
\end{equation}
Let us also notice that from (\ref{eq:mix_time}) it immediately follows that, for every $C>0$, there exists $\kappa$ such that, for all $t\geq\kappa\ln N$ and all $x,y\in[N]$, we have 
\begin{equation}
\left|P^t(x,y)-\frac{1}{N}\right|<\frac{1}{N^2}.
\label{eq:mixing_Janson}
\end{equation}

Let us fix a positive $\tilde c<\min\{c,1\}$ and prove that the random walk subject to $X_t=X_0$ intersects itself during the first $\tilde{c}\sqrt{N}$ steps with non-zero probability.
Let
$$
\mathcal{J}:=\setdef{(i,j)\in[t]^2}{\kappa\ln N<i<j-3\kappa\ln N<j<\tilde{c}\sqrt{N}}.
$$
Note that $|\mathcal{J}|=\frac{\tilde c^2+o(1)}{2}N$
and that $j<t-\kappa\ln{n}$ for every $(i,j)\in\mathcal{J}$.

For all $(i,j)\in\mathcal{J}$, we have that
$$
\frac{\Pb{X_i=X_j,X_t=X_0}}{\Pb{X_t=X_0}}=
\sum_{x,u\in [N]} \frac{P^i(x,u)P^{j-i}(u,u)P^{t-j}(u,x)}{\sum_{y\in [N]}P^t(y,y)}
=
\frac{1+o(1)}{N}
$$
due to (\ref{eq:mixing_Janson}).
Note that $o(1)$ in the expression above converges to $0$ uniformly over all $(i,j)\in\mathcal{J}$.
Therefore,
\begin{equation}
\sum_{(i,j)\in\mathcal{J}}\frac{\Pb{X_i=X_j,X_t=X_0}}{\Pb{X_t=X_0}}
=
\frac{\tilde c^2}{2}+o(1).
\label{eq:relative_error_first_moment}
\end{equation}
Let
$$
\begin{aligned}
\mathcal{J}_0&:=
\setdef{\{(i_1,j_1),(i_2,j_2)\}\in{\mathcal{J}\choose 2}}{\pr{\forall\br{i,j}\subset\{i_1,i_2,j_1,j_2\}:\abs{i-j}\ge\kappa\ln{N}}},\\
\mathcal{J}_1&:=
\setdef{\{(i_1,j_1),(i_2,j_2)\}\in{\mathcal{J}\choose 2}}{\pr{\exists!\br{i,j}\subset\{i_1,i_2,j_1,j_2\}:\abs{i-j}<\kappa\ln{N}}},\\
\mathcal{J}_2&:=
\setdef{\{(i_1,j_1),(i_2,j_2)\}\in{\mathcal{J}\choose 2}}{\max\br{\abs{i_1-i_2},\abs{j_1-j_2}}<\kappa\ln{N}}.
\end{aligned}
$$
It is clear from the definition of $\mathcal{J}$ that $\mathcal{J}_0\sqcup\mathcal{J}_1\sqcup\mathcal{J}_2={\mathcal{J}\choose 2}$ 
(recall that
$j-i>3\kappa\ln{N}$ for every $(i,j)\in\mathcal{J}$).
We have 
$$
|\mathcal{J}_0|=\frac{1}{2}|\mathcal{J}^2|(1+o(1))=\frac{\tilde c^4+o(1)}{8}N^2,\quad
|\mathcal{J}_1|<4\kappa\ln N(\tilde c\sqrt{N})^3,\quad
|\mathcal{J}_2|<(2\kappa\ln N\tilde c\sqrt{N})^2.
$$
As above, uniformly over all $\{(i_1,j_1),(i_2,j_2)\}\in\mathcal{J}_0$, \eqref{eq:mixing_Janson} implies
$$
\frac{\Pb{X_{i_1}=X_{j_1},X_{i_2}=X_{j_2},X_t=X_0}}{\Pb{X_t=X_0}}
=
\frac{1+o(1)}{N^2}.
$$
Uniformly over all $\{(i_1,j_1),(i_2,j_2)\}\in\mathcal{J}_1$, the relations (\ref{eq:transition_above}) and (\ref{eq:mixing_Janson}) imply
$$
\frac{\Pb{X_{i_1}=X_{j_1},X_{i_2}=X_{j_2},X_t=X_0}}{\Pb{X_t=X_0}}
\leq\frac{1+o(1)}{NN_1}.
$$
Uniformly over all $\{(i_1,j_1),(i_2,j_2)\}\in\mathcal{J}_2$, \eqref{eq:transition_above} implies
$$
\frac{\Pb{X_{i_1}=X_{j_1},X_{i_2}=X_{j_2},X_t=X_0}}{\Pb{X_t=X_0}}
\leq\frac{1+o(1)}{N_1^2}.
$$
Summing up and recalling that, in the current case, $\sqrt{N}=o(N_1)$,
\begin{multline}
\sum_{\{(i_1,j_1),(i_2,j_2)\}\in{\mathcal{J}\choose 2}}
\frac{\Pb{X_{i_1}=X_{j_1},X_{i_2}=X_{j_2},X_t=X_0}}{\Pb{X_t=X_0}}\\
\leq \frac{\tilde c^4+o(1)}{8}+\frac{(4\kappa\tilde c^3+o(1))\sqrt{N}\ln N}{N_1}+
\frac{(2\kappa\tilde c+o(1))^2 N \ln N}{N_1^2}=
\frac{\tilde c^4+o(1)}{8}.
\label{eq:relative_error_second_moment}
\end{multline}
From (\ref{eq:relative_error}), (\ref{eq:relative_error_first_moment}) and (\ref{eq:relative_error_second_moment}) we get
\begin{align*}
\frac{\homo{C_t,G}-\mono{C_t,G}}{\homo{C_t,G}} & \geq
\sum_{(i,j)\in\mathcal{J}}
\frac{\Pb{X_i=X_j,X_t=X_0}}{\Pb{X_t=X_0}}\\
&-\sum_{\{(i_1,j_1),(i_2,j_2)\}\in{\mathcal{J}\choose 2}}
\frac{\Pb{X_{i_1}=X_{j_1},X_{i_2}=X_{j_2},X_t=X_0}}{\Pb{X_t=X_0}}\\
&=\frac{\tilde c^2}{2}-\frac{\tilde c^4}{8}+o(1).
\end{align*}
Since $\frac{\tilde c^2}{2}-\frac{\tilde c^4}{8}>0$, we conclude that $\mono{C_t,G}/\homo{C_t,G}$ is bounded away from $1$ as needed.

\subsection{$N_1=O(\sqrt{N})$}
\label{sec:th3_proof_2}
W.l.o.g. we may assume that $c<1$ and prove that the random walk subject to $X_t=X_0$ intersects itself during the first $cN_1$ steps with non-zero probability. In the same way as in Section~\ref{sec:th3_proof_1}, we use~(\ref{eq:relative_error}). However, here we consider all $(i,j)$ such that $i$ is even and $j=i+2\leq cN_1$. Let $\mathcal{J}:=\setdef{2i}{i\in\irange{0}{\floor{\frac{cN_1-2}{2}}}}$. For every
$i\in\mathcal{J}
$,
we have
$$
\frac{\Pb{X_i=X_{i+2},X_t=X_0}}{\Pb{X_t=X_0}}
=
\frac{1}{N}\sum_{x,u\in[N]} \frac{P^i(x,u)P^2(u,u)P^{t-i-2}(u,x)}{\Pb{X_t=X_0}}
=
\frac{1}{N_1}\frac{\Pb{X_{t-2}=X_0}}{\Pb{X_t=X_0}}
\sim
\frac{1}{N_1}
$$
since $P^2(u,u)=\frac{1}{N_1}$ for all $u\in[N]$ and $\Pb{X_{t-2}=X_0}\sim\Pb{X_t=X_0}\sim1/N$ due to (\ref{eq:mixing_Janson}). Thus,
$$
\sum_{i\in\mathcal{J}}
\frac{\Pb{X_i=X_{i+2},X_t=X_0}}{\Pb{X_t=X_0}}
=
\frac{c}{2}+o(1).
$$
Now, let $i_1,i_2\in\mathcal{J}$, $i_1<i_2$. Then, similarly,


\begin{align*}
&\frac{\Pb{X_{i_1}=X_{i_1+2},X_{i_2}=X_{i_2+2},X_t=X_0}}{\Pb{X_t=X_0}}=\\
&=\frac{1}{N}\sum_{x,u,v\in[N]} \frac{P^{i_1}(x,u)P^2(u,u)P^{i_2-i_1-2}(u,v)P^2(v,v)P^{t-i_2-2}(v,x)}{\Pb{X_t=X_0}}\\
&=\frac{1}{N}\sum_{x,u,v\in[N]}
\frac{P^{i_1}(x,u)P^{i_2-i_1-2}(u,v)P^{t-i_2-2}(v,x)}{N_1^2\,\Pb{X_t=X_0}}\\
&=\frac{\Pb{X_{t-4}=X_0}}{N_1^2\,\Pb{X_t=X_0}}
\sim\frac{1}{N_1^2}.
\end{align*}

Therefore,
\begin{align*}
\frac{\homo{C_t,G}-\mono{C_t,G}}{\homo{C_t,G}} & \geq
\sum_{i\in\mathcal{J}}
\frac{\Pb{X_i=X_{i+2},X_t=X_0}}{\Pb{X_t=X_0}}\\
&-\sum\limits_{i_1,i_2\in\mathcal{J}:i_1<i_2}
\frac{\Pb{X_{i_1}=X_{i_1+2},X_{i_2}=X_{i_2+2},X_t=X_0}}{\Pb{X_t=X_0}}\\
&=\frac{c}{2}-\frac{c^2}{8}+o(1).
\end{align*}
Since $c/2-c^2/8>0$, this finishes the proof of Theorem~\ref{th_mon_asymp_non_eq_hom} and therefore of Theorem~\ref{th_mon_asymp_eq_hom}. \qed

\section{Acknowledgements}

We would like to thank Andrey Kupavskii for very useful discussions.\\


\begin{thebibliography}{99}


\bibitem{Agong2018} L. A. Agong, C. Amarra, J. S. Caughman, A. J. Herman, T. S. Terada, {\it On the girth and diameter of generalized Johnson graphs}, Discrete Mathematics, {\bf 341}:1 (2018) pp. 138--142.

\bibitem{Alspach2012} B. Alspach, {\it Johnson graphs are Hamilton-connected}, Ars Mathematica Contemporanea, {\bf 6}:1 (2012) pp. 21--23.

\bibitem{Bapat2014} R. B. Bapat, \textit{Graphs and Matrices}, Springer London (2014).

\bibitem{Brouwer1979} A. E. Brouwer, A. Schrijver, {\it Uniform Hypergraphs}, In  Packing and Covering in Combinatorics. Mathematical Centre Tracts, {\bf 106} (1979) pp. 39--73.

\bibitem{Burkin2016} A. V. Burkin, {\it Small Subgraphs in Random Distance Graphs}, Theory of Probability {\&} Its Applications, {\bf 60}:3 (2016) pp. 367--382.

\bibitem{Burkin2018} A. V. Burkin, M. E. Zhukovskii, {\it Small subgraphs and their extensions in a random distance graph},  Sbornik: Mathematics, {\bf 209}:2 (2018) pp. 163--186.

\bibitem{Cannon2012} A. D. Cannon, J. Bamberg, C. E. Praeger, {\it A Classification of the Strongly Regular Generalised Johnson Graphs}, Annals of Combinatorics, {\bf 16}:3 (2012), pp.489--506.

\bibitem{Cantwell1996} K. Cantwell, {\it Finite Euclidean Ramsey theory}, Journal of Combinatorial Theory,  Series A, {\bf 73}:2 (1996) pp. 273--285.


\bibitem{Chen1987} B.-L. Chen, K.-W. Lih, {\it Hamiltonian uniform subset graphs}, Journal of Combinatorial Theory,  Series B, {\bf 42}:3 (1987) 257--263.

\bibitem{Chen2003} Y.-C. Chen, {\it Triangle-free Hamiltonian Kneser graphs}, Journal of Combinatorial Theory,  Series B, {\bf 89}:1 (2003) pp. 1--16.

\bibitem{Chen2008} Y. Chen, W. Wang, {\it Diameters of uniform subset graphs}, Discrete Mathematics, {\bf 308}:24 (2008) 6645--6649.

\bibitem{Chen2008_generalized_kneser} Y. Chen, Y. Wang, {\it On the diameter of generalized Kneser graphs},  Discrete Mathematics, {\bf 308}:18 (2008) 4276--4279.

\bibitem{Chilakamarri1990} K. B. Chilakamarri, {\it On the chromatic number of rational five-space}, Aequationes Mathematicae, {\bf 39}:2-3 (1990) pp. 146--148.

\bibitem{Daven1999} M. Daven, C. A. Rodger, {\it The Johnson graph $J(v,k)$ has connectivity $\delta$}, Congressus Numerantium, {\bf 139} (1999) pp. 123--128.

\bibitem{Delsarte1973} P. Delsarte, \textit{An algebraic approach to the association schemes of coding theory}, Philips Res. Rep. Suppl. No. 10 (1973).

\bibitem{Denley1997} T. Denley, {\it The Odd Girth of the Generalised Kneser Graph}, European Journal of Combinatorics, {\bf 18}:6 (1997) pp. 607--611.

\bibitem{Etzion1996_codes} T. Etzion, {\it On the Nonexistence of Perfect Codes in the Johnson Scheme}, 
{SIAM} Journal on Discrete Mathematics, {\bf 9}:2 (1996) pp. 201--209.

\bibitem{Etzion1996_chromatic} T. Etzion, S. Bitan, {\it On the chromatic number, colorings, and codes of the Johnson graph}, Discrete Applied Mathematics, {\bf 70}:2 (1996) pp. 163--175.

\bibitem{Etzion2011} T. Etzion, A. Brouwer, {\it Some new distance-4 constant weight codes}, Advances in Mathematics of Communications, {\bf 5}:3 (2011) 417--424.

\bibitem{Exoo2014} G. Exoo, D. Ismailescu, M. Lim, {\it On the Chromatic Number of $\mathbb{R}^4$}, Discrete {\&} Computational Geometry, {\bf 52}:2 (2014) pp. 416--423.

\bibitem{Frankl1980} P. Frankl, {\it Extremal Problems and Coverings of the Space}, European Journal of Combinatorics, {\bf 1}:2 (1980) pp. 101--106.


\bibitem{Frankl1985_generalized_kneser} P. Frankl, {\it On the chromatic number of the general Kneser-graph}, Journal of Graph Theory, {\bf 9}:2 (1985) pp. 217--220.



\bibitem{Frankl1981} P. Frankl, R. M. Wilson, {\it Intersection theorems with geometric consequences}, Combinatorica, {\bf 1}:4 (1981) pp. 357--368.


\bibitem{Gao2021} P. Gao, \textit{Kim-Vu's sandwich conjecture is true for all $d=\omega(\log^7 n)$}, arXiv:2011.09449 (2021).

\bibitem{Jafari2020} A. Jafari, M. J. Moghaddamzadeh, {\it On the chromatic number of generalized Kneser graphs and Hadamard matrices}, Discrete Mathematics, {\bf 343}:2 (2020) 111682.

\bibitem{Janson2000} S.~Janson, T.~Luczak, A.~Ruci{\'{n}}ski, \textit{Random Graphs}, Wiley, New York (2000).

\bibitem{Kahn1993} J. Kahn, G. Kalai, {\it A Counterexample to Borsuk's Conjecture}, Bulletin of the American Mathematical Society, {\bf 29}:1 (1993) pp. 60--63.

\bibitem{Kupavskii2009} A. B. Kupavskii and A. M. Raigorodskii, {\it On the chromatic number of $\mathbb{R}^9$}, Journal of Mathematical Sciences, {\bf 163}:6 (2009) 720--731.

\bibitem{Larman1972} D. G. Larman, C. A. Rogers, {\it The realization of distances within sets in Euclidean space}, Mathematika, {\bf 19}:1 (1972) pp. 1--24.


\bibitem{Lovasz1978} L. Lov{\'{a}}sz, {\it Kneser's conjecture, chromatic number, and homotopy}, Journal of Combinatorial Theory,  Series A, {\bf 25}:3 (1978)
 pp. 319--324.
 
\bibitem{Matousek2004} J. Matou{\v{s}}ek, {\it A Combinatorial Proof of Kneser's Conjecture}, Combinatorica, {\bf 24}:1 (2004) pp. 163--170.
 
\bibitem{Meyer2000} C. D. Meyer, \textit{Matrix analysis and applied linear algebra}, Society for Industrial and Applied Mathematics, Philadelphia (2000).
 
\bibitem{Molitierno2017} J. J. Molitierno, {\it Entries of the group inverse of the Laplacian matrix for generalized Johnson graphs}, Linear and Multilinear Algebra, {\bf 66}:6 (2017) pp. 1153--1170.

\bibitem{Mutze2020} T. M\"{u}tze, J. Nummenpalo, B. Walczak, {\it Sparse Kneser graphs are Hamiltonian}, Journal of the London Mathematical Society {\bf 103}:4 (2021) 1253--1275.


\bibitem{Poljak1987} S. Poljak, Z. Tuza, {\it Maximum bipartite subgraphs of Kneser graphs}, Graphs and Combinatorics, {\bf 3}:1 (1987) pp. 191--199.


\bibitem{Pyaderkin2016} M. M. Pyaderkin, {\it Independence numbers of random subgraphs of distance graphs}, Mathematical Notes, {\bf 99}:3-4 (2016) pp. 556--563.


\bibitem{Raigorodskii2001} A. M. Raigorodskii, {\it Borsuk's problem and the chromatic numbers of some metric spaces}, Russian Mathematical Surveys, {\bf 56}:1 (2001) pp. 103--139.


\bibitem{Raigorodskii2016} A. M. Raigorodskii, {\it Combinatorial Geometry and Coding Theory}, Fundamenta Informaticae, {\bf 145}:3 (2016) pp. 359--369.

\bibitem{Raigorodskii2013} A. M. Raigorodskii, {\it Coloring Distance Graphs and Graphs of Diameters}, Thirty Essays on Geometric Graph Theory ed. J{\'{a}}nos Pach (2013) pp. 429--460.

\bibitem{Shields2002} I. Shields, C. Savage, \textit{A Note on Hamilton Cycles in Kneser Graphs}, Bull. Inst. Combin. Appl. \textbf{40} (2002).

\bibitem{Simpson1994} J. E. Simpson, {\it On uniform subset graphs}, Ars Combinatoria, {\bf 37} (1994) pp. 309--318.

\bibitem{Sinclair1992} A. Sinclair, \textit{Improved bounds for mixing rates of Markov chains and multicommodity flow}, in: LATIN 1992. Lecture Notes in Computer Science, \textbf{583} (1992).

\bibitem{Valencia2005} M. Valencia-Pabon, J.-C. Vera, {\it On the diameter of Kneser graphs}, Discrete Mathematics, {\bf 305}:1-3 (2005) pp. 383--385.

\bibitem{Zhao2012} Y. Zhao, \textit{Spectral Distributions of Random Graphs} (2012).

\bibitem{Zhukovskii2012_sub} M. E. Zhukovskii, {\it On the probability of the occurrence of a copy of a fixed graph in a random distance graph}, Mathematical Notes, {\bf 92}:6 (2012) pp. 756--766.



\end{thebibliography}

\end{document}